\documentclass[reqno,12pt]{amsart}
\usepackage{amsmath, amssymb, amsthm, amsfonts} 
\usepackage[english]{babel}
\usepackage{bbm}
\usepackage{graphicx}
\usepackage{url}
\usepackage{soul}
\usepackage{epstopdf}
\usepackage{float}
\usepackage[ruled,vlined]{algorithm2e}
\usepackage[a4paper,bindingoffset=0.5cm,left=2cm,right=2cm,top=2.5cm,bottom=2cm,footskip=.8cm]{geometry}
\usepackage{rotating}
\usepackage{amsbsy,enumerate}
\usepackage{comment}
\usepackage{mathrsfs} 
\usepackage{amsmath, amssymb, amsthm,amsfonts}

\newcommand{\LMG}{\xi}
\newcommand{\LM}{\lambda+\mu}
\newcommand{\MG}{\mu+\gamma}
\newcommand{\LG}{\lambda+\gamma}

\renewcommand{\hat}{\widehat}

\newcommand{\vertiii}[1]{{\left\vert\kern-0.25ex\left\vert\kern-0.25ex\left\vert #1 \right\vert\kern-0.25ex\right\vert\kern-0.25ex\right\vert}}
\allowdisplaybreaks
\usepackage{hyperref}

\allowdisplaybreaks

\newtheorem{theorem}{Theorem}
\newtheorem{proposition}{Proposition}

\theoremstyle{remark}
\newtheorem{remark}{Remark}

\usepackage[utf8]{inputenc}
\usepackage{type1cm}         
\usepackage{multicol}        
\usepackage[bottom]{footmisc}

\usepackage{newtxtext}       %
\usepackage{newtxmath}       

\usepackage{tikz}
\usetikzlibrary{plotmarks}
\usetikzlibrary{decorations.pathreplacing}
\usepackage{pgfplots}
\pgfplotsset{width=10cm,compat=1.9}
\usepgfplotslibrary{external}
\usepackage{subcaption}
\usepackage{standalone}
\usetikzlibrary{automata,positioning}
\usetikzlibrary{decorations.pathreplacing,decorations.markings,shapes.geometric}
\usetikzlibrary{shapes,arrows}
\usetikzlibrary{backgrounds,calc,positioning}

\usepackage{ragged2e}
\usepackage{booktabs,tabularx}

\usepackage{pgfplots}
\pgfplotsset{compat=1.9}

\usetikzlibrary{chains,shapes.multipart}
\usetikzlibrary{shapes,calc}
\usetikzlibrary{automata,positioning}

\usepackage{tikz}
\usetikzlibrary{automata,positioning}
\usetikzlibrary{decorations.pathreplacing,decorations.markings,shapes.geometric}
\usetikzlibrary{shapes,arrows}
\usetikzlibrary{backgrounds,calc,positioning}

\newcommand{\eq}[1]{(\ref{eq:#1})}

\begin{document}

	\title[Finite customer-pool queues]{Finite customer-pool queues}
\author{Onno Boxma, Offer Kella,
{\tiny and} Michel Mandjes}
	
    \begin{abstract}
    In this paper we consider an M/G/1-type queue fed by a finite customer-pool. 
    In terms of transforms, we characterize the time-dependent distribution of the number of customers and the workload, as well as the associated waiting times.

\smallskip


\smallskip

\noindent
{\sc Affiliations.} 
OB is with {E}{\sc urandom} and Department of Mathematics and Computer Science, Eindhoven University of Technology, the Netherlands.
OK is with Department of Statistics and Data Science, The Hebrew University of Jerusalem; Jerusalem 9190501, Israel. MM is with the Mathematical Institute, Leiden University, 
The Netherlands. He is also affiliated with Korteweg-de Vries Institute for Mathematics, University of Amsterdam, The Netherlands; E{\sc urandom}, Eindhoven University of Technology, The Netherlands; Amsterdam Business School, University of Amsterdam, The Netherlands.


\noindent

\smallskip

\noindent
{\sc Acknowledgments.} 
The research of OB and MM was supported by the European Union’s Horizon 2020 research and innovation programme under the Marie Sklodowska-Curie grant agreement no.\ 945045, and by the NWO Gravitation project NETWORKS under grant agreement no.\ 024.002.003. The research of OK is supported in part by the NSF grant~3336/24 and by the Vigevani Chair in Statistics.
Date: {\it \today}.

\smallskip

\noindent
{\sc Email.} {\scriptsize \url{o.j.boxma@tue.nl}, \url{offer.kella@gmail.com}, \url{m.r.h.mandjes@math.leidenuniv.nl}}

	\end{abstract}

	\maketitle

\section{Introduction}
The model that we consider in this paper can be interpreted as a finite customer-pool M/G/1 system. It has the distinguishing feature that, after time 0, only a finite number $m$ of customers bring work into the system, while there are an additional $k$ customers present at time $0$ (whose service has not started at time $0$). Once a customer's work has been processed, it leaves the system, entailing that the queue eventually empties with probability 1 (i.e., without any stability requirement being imposed). The model is of M/G/1 type in that in our setup the interarrival times are independent and exponentially distributed with a parameter that depends on the number of customers yet to arrive. 

The objective of this study is to describe, for the model introduced in the previous paragraph, the transient behavior of the number of customers present and the workload. Indeed, we succeed in finding a recursive algorithm (of low complexity) that provides us with the joint transform of these two random quantities at an exponentially distributed time. Our findings allow the evaluation of the joint probability mass function (for the number of customers) and density (for the workload) at any given point in time, relying on e.g.\ the well-developed machinery of \cite{AW,dI}. 

Our work directly relates to two branches of work. (i)~In the first place, there is a connection to papers in which the {\em arrival} times (rather than {\em interarrival} times) are independent, identically distributed (i.i.d.) random variables. Indeed, when taking in our model the parameter of the exponentially distributed interarrival times proportional to the number of customers yet to arrive, we can equivalently say that the arrival times are i.i.d.\ exponentially distributed random variables. 
This setup with independent arrival times is particularly appropriate in instances that the number of potential customers is relatively low, so that the usual Poisson arrival assumption cannot be justified.  
The resulting model has been analyzed in \cite{BET,BHL,HON}, with emphasis on scaling limits in the diffusion regime; we in addition refer to the survey \cite{RH}. The recent study \cite{RUT} provides an exact analysis of the workload at an exponentially distributed time, also covering various model variants (additional Poisson arrival stream, balking customers, etc.). (ii)~Our model can also be seen as a special case of a reflected Markov additive process \cite[Ch.\ XI]{ASM2} in which the underlying background process is transient. More concretely, the state of the background process represents the number of customers yet to arrive, a process that decreases over time. The queueing model driven by a MAP with a transient background process has been analyzed in \cite{KMD}; see also  \cite{DelMan} for its ruin-theoretic counterpart. In the present paper we show for our specific model, in which the background process is of pure-death type, that the expressions of \cite{KMD} simplify considerably. 

Queues with independent arrival times form a natural class of models in situations where customers independently decide when to try to access service. The arrival time distribution could e.g.\ reflect the underlying {\em day profile}. As pointed out in e.g.\ \cite{RH}, models of this type lead to intriguing game-theoretic questions; this stream of the literature goes back to the seminal work \cite{NAOR}. In our paper we explicitly include the special situation in which at time 0 already some jobs, characterized via their specific service time distribution, can be present; this model instance can be used if the server has to perform work prior to being able to serve the customers who arrive after time 0.

\smallskip

We conclude this introduction by detailing our contributions and the organization of the paper. In Section \ref{model} we present the model and define the key objects of interest. Our first main result, in Section \ref{sec:ana}, concerns the distribution of the number of customers at an exponentially distributed time (which, by applying inversion, gives us the distribution at any deterministic time). It is also pointed out how the result can be extended to the joint distribution with the workload. Our findings refute the claim in \cite{HON} that the finite customer-pool M/G/1 system be `not amenable to exact analysis', but the obvious price we pay is that we have to settle for results in the transform domain.  Then Section \ref{sec:wt} presents an analysis of the waiting-time distribution, including tail asymptotics in case the customers' service times are regularly varying. Section \ref{sec:spec} presents several ramifications. 


\section{Model, notation, and objective}\label{model}
In this section we formally introduce our model, state our objectives, and define two sequences of probabilities that play a crucial role throughout our analysis. 
We throughout denote $a^+=\max(a,0)$, $a^-=\min(a,0)$, LST abbreviates {\em Laplace Stieltjes Tranform}, `a.s.' abbreviates {\em almost surely}, $B$ is some a.s.\ positive (and finite) random variable with LST $\beta(\alpha)=\mathbb{E}\,{\rm e}^{-\alpha B}$ for $\alpha\geqslant 0$. We follow the convention that ${\mathbb N}$ and $\mathbb{N}_0$ denote the sets of, respectively, positive and nonnegative integers.

In our model $k\in{\mathbb N}_0$ customers are already present at time $0$, but their service times have not started at time $0$. In addition, $m\in{\mathbb N}_0$ customers are to arrive after time $0$. With ${\rm Exp}(\lambda)$ for $\lambda>0$ denoting an exponentially distributed random variable with mean $\lambda^{-1}$, the random variable $T_m\sim{\rm Exp}(\lambda_m)$ is the time until the first arrival, $T_{m-1}\sim{\rm Exp}(\lambda_{m-1})$ is the time between the first and second arrivals, etc. 
All $k+m$ service times are independent (and independent of the arrival process), and distributed like $B$. Observe that this model can be seen as a finite customer-pool M/G/1 queue.

In the sequel, ${\mathbb E}_{k,m}[\cdot]$ denotes the expectation when starting with $k\in{\mathbb N}_0$ customers being present, and $m\in{\mathbb N}_0$ customers to arrive according to the mechanism described above. The main quantity of interest is $Z(t)$, the number of customers present at time $t\geqslant 0$. Let $T\sim{\rm Exp}(\gamma )$, independently of everything else. It is our objective to compute the generating function \begin{equation}
\mu_{km}(z)={\mathbb E}_{k,m}[z^{Z(T)}],\quad z\in (0,1],  
\end{equation} of the number of customers present at the exponentially distributed time $T$. Evidently,
\begin{equation}
    \mu_{km}(z)=\int_0^\infty \gamma  {\rm e}^{-\gamma  t}\,{\mathbb E}_{k,m}[z^{Z(t)}] \,{\rm d}t,
\end{equation}
so that one can numerically obtain ${\mathbb E}_{k,m}[z^{Z(t)}]$ by applying Laplace inversion; see for instance the fast and accurate procedures presented in  \cite{AW,dI}. By repeated differentiation and inserting $0$, one can obtain the probability distribution of $Z(t)$, while repeated differentiation and inserting 1 yields all (factorial) moments. 

\medskip

In our analysis two sequences of auxiliary probabilities play a key role. 
To this end, we denote $S_0=0$ and $S_n:=\sum_{i=1}^n T_i$ for any $n\in\{1,\ldots,m\}$.
\begin{itemize}
    \item[$\circ$]
The first sequence of probabilities is defined as follows. Let $u_{00}:={\mathbb P}(T>B)=\beta(\gamma)$ and for $n\in\{0,\ldots,m\}$,
\begin{align}\label{eq:umi}
u_{n i}:&=
\begin{cases}
{\mathbb P}(S_n-S_{n-i}\leqslant B<S_n-S_{n-i-1},T>B),&i\in\{0,\ldots, n-1\}\,,\\
{\mathbb P}(S_n\leqslant B,T>B),&i=n\,,
\end{cases}\nonumber\\
&=
\begin{cases}
{\mathbb E}\left({\rm e}^{-\gamma B}\,1_{\{S_n-S_{n-i}\leqslant B<S_n-S_{n-i-1}\}}\right),&\hskip 1.4 cm i\in\{0,\ldots, n-1\}\,,\\
{\mathbb E}\left({\rm e}^{-\gamma B}\,1_{\{S_n\leqslant B\}}\right),&\hskip 1.4 cm i=n\,,
\end{cases}
\end{align}
noting that $\sum_{i=0}^nu_{ni}=\mathbb{P}(T>B)=\beta(\gamma)$.
The object $u_{ni}$ can be interpreted as the probability that there are $i$ arrivals during a single service time, starting with $n\in\{0,\ldots,m\}$ customers who are yet to arrive, jointly with the event $\{T>B\}$.
 \item[$\circ$]
The second sequence of probabilities is defined similarly: we let $v_{00}:={\mathbb P}(T\leqslant B)=1-\beta(\gamma)$, while for $n\in\{0,\ldots,m\}$,
\begin{equation}\label{eq:vmi}
v_{n i}:=
\begin{cases}
{\mathbb P}(S_n-S_{n-i}\leqslant T<S_n-S_{n-i-1},T\leqslant B),&i\in\{0,\ldots, n-1\},\\
{\mathbb P}(S_n\leqslant T,T\leqslant B ),&i=n\,,
\end{cases}
\end{equation}
noting that $\sum_{i=0}^nv_{ni}=\mathbb{P}(T\leqslant B)=1-\beta(\gamma)$.
In this case, $v_{n i}$ can be interpreted as the probability that there are $i$ arrivals until the exponential time $T$, jointly with $T\leqslant B$.
{
We observe that for $i\in\{0,\ldots,n-1\}$ we can write
\begin{equation}\label{eq:vmi1}
v_{ni}={\mathbb P}(S_n-S_{n-i}\leqslant T<S_n-S_{n-i-1})-{\mathbb P}(S_n-S_{n-i}\leqslant T<S_n-S_{n-i-1},T> B)
\end{equation}
where
\begin{equation}\label{eq:vmi2}
    {\mathbb P}(S_n-S_{n-i}\leqslant T<S_n-S_{n-i-1})=\frac{\gamma}{\lambda_{n-i}+\gamma}\prod_{j=n-i+1}^n\frac{\lambda_j}{\lambda_j+\gamma}
\end{equation}
(with an empty product defined to be~$1$) and, by the memoryless property,
\begin{align}\label{eq:vmi3}
    {\mathbb P}(S_n-S_{n-i}&\leqslant T<S_n-S_{n-i-1},T> B)
    =\,{\mathbb E}[{\rm e}^{-\gamma B}1_{\{S_n-S_{n-i}\leqslant B+T<S_n-S_{n-i-1}\}}]\,.
\end{align}
Along the same lines,
\begin{equation}\label{eq:vmi4}
    v_{nn}=\prod_{j=1}^n\frac{\lambda_j}{\lambda_j+\gamma}-\mathbb{E} [{\rm e}^{-\gamma B}1_{\{S_n\leqslant B+T\}}]\,.
\end{equation}
}
\end{itemize}
\begin{remark}\label{remBexp}
We note that when $B\sim \text{Exp}(\mu)$ then $B\wedge T$ and $\{T\leqslant B\}$ (as well as $B\wedge T$ and $\{B>T\}$) are independent, with $B\wedge T\sim \text{Exp}(\gamma+\mu)$ and ${\mathbb P}(B\leqslant  T)=1-{\mathbb P}(B>T)=\frac{\mu}{\gamma+\mu}$. Therefore,
\begin{align}
\frac{u_{ni}}{\mu}=\frac{v_{ni}}{\gamma}&=\frac{1}{\gamma+\mu}\,\mathbb{P}\left(S_n-S_{n-i}\leqslant  B\wedge T<S_n-S_{n-i-1}\right)\nonumber \\
&=\frac{1}{\gamma+\mu}\,\frac{\gamma+\mu}{\lambda_{n-i}+\gamma+\mu}\,
\prod_{j=n-i+1}^n\frac{\lambda_j}{\lambda_j+\gamma+\mu}\\
&=\frac{1}{\lambda_{n-i}+\gamma+\mu}\,
\prod_{j=n-i+1}^n\frac{\lambda_j}{\lambda_j+\gamma+\mu},
\end{align}
defining $\lambda_0=0$.
$\hfill\Diamond$
\end{remark}
In Section  \ref{sec:ana} the probabilities $u_{ni}$ and $v_{ni}$ are, in principle, considered known. Observe in particular that they are functions of the parameter $\gamma $ corresponding to the `killing time' $T$. In Section \ref{sec:spec} we consider special instances in which  $u_{ni}$ and $v_{ni}$ allow closed-form expressions.

\section{Number of customers at killing} \label{sec:ana}
In the analysis of the conventional M/G/1 queue, one typically works with the embedded process that records the number of customers present at departure epochs. This process is a Markov chain: if it has value $\ell\in {\mathbb N}$, then at the next epoch it is $\ell-1$ increased by the number of Poisson arrivals during the service time of the served customer.

In our model we exploit a similar idea. 
Observe that if the number of customers at a certain departure epoch is $\ell\in{\mathbb N}$ and $T>B$, then the number of customers at the next departure epoch will be given by $\ell-1$ plus the number of arrivals during the corresponding service time (say $i$). In addition, supposing that there were still $n\in\{0,\ldots,m\}$ customers yet to arrive, this number becomes $n-i$. Therefore,  we have from the memoryless property for $T$ and the customers that have not arrived until time $B$, for any $\ell\in\{1,\ldots,k\}$ that
\begin{equation}\label{C1}
{\mathbb E}_{\ell,n}[z^{Z(T)}1_{\{T>B\}}]=\sum_{i=0}^n \mu_{\ell+i-1,n-i}(z)\,u_{n,i}\,.
\end{equation}
We proceed by distinguishing various cases. As we point out, the resulting expressions allow a recursive algorithm by which $\mu_{km}(z)$ can be evaluated. 
\begin{itemize}
    \item[$\circ$] 
    When $\ell\in\mathbb{N}$ and $T\leqslant B$, then the total number of customers in the system at time $T$ is $\ell$ plus the number of customers that have arrived by time $T$. Therefore in this case, for any $n\in\{0,\ldots,m\}$,
\begin{equation}\label{C2}
{\mathbb E}_{\ell,n}[z^{Z(T)}1_{\{T\leqslant B\}}]=z^\ell\sum_{i=0}^n z^i\,v_{ni}\,.
\end{equation}
Upon combining the displays \eqref{C1} and \eqref{C2},  this gives the relation, for any $\ell\in\{ 1,\ldots,k\}$ and $n\in\{0,\ldots,m\}$,
\begin{equation}\label{eq:mukm}
\mu_{\ell n}(z)=\sum_{i=0}^n \left(\mu_{\ell+i-1,n-i}(z)\,u_{ni}+ z^{\ell+i}\,v_{ni}\right)\,.
\end{equation}
    \item[$\circ$] 
    Supposing that $\ell=0$, $n\in\{1,\ldots,m\}$  and $T\leqslant T_n$, we evidently have that $Z(T)=0$. In case $\ell=0$, $n\in\{1,\ldots,m\}$  and $T>T_n$, then by the memoryless property, $Z(T)$ behaves as if we would have started from the state $(1,n-1)$. As a result,
\begin{equation}\label{eq:zerom}
\mu_{0n}(z)=\frac{\gamma }{\gamma +\lambda_n}+\frac{\lambda_n}{\gamma +\lambda_n}\,\mu_{1,n-1}(z)\,.
\end{equation}
    \item[$\circ$] Finally, when starting with $\ell=n=0$, then clearly $Z(T)=0$ and thus $\mu_{00}(z)=1$. For $n=0$ and $\ell\in\{ 1,\ldots,k\}$, \eq{mukm} becomes
\begin{equation}\label{eq:kzero0}
\mu_{\ell 0}(z)=\mu_{\ell-1,0}(z)\,u_{00}+z^\ell \,v_{00}\,.
\end{equation}
This implies, by induction, that
\begin{equation}\label{eq:kzero}
\mu_{\ell 0}(z)=
\begin{cases}
(1-v_{00})u_{00}^\ell +v_{00}\,{\displaystyle \frac{u_{00}^{\ell +1}-z^{\ell +1}}{u_{00}-z}},&z\not=u_{00},\\
(1+v_{00}\ell)\,u_{00}^\ell ,&z=u_{00}\,.
\end{cases}
\end{equation}
\end{itemize}

The above can be used to devise a procedure that recursively determines $\mu_{km}(z)$. It is given by the following algorithm.

\medskip
\begin{quotation} {\footnotesize
{\tt ALGORITHM to evaluate $\mu_{km}(z)$}.

\noindent {\tt INPUT: $k$, $m$, $u_{ni}$ and $v_{ni}$ for all $n\in
\{0,\ldots,m\}$ and $i\in\{0,\ldots,n\}$.}

\noindent {\tt OUTPUT: $\mu_{\ell n}(z)$ for all $\ell\in\{0,\ldots,k\}$ and $n\in\{0,\ldots,m\}.$}

\medskip

    \noindent {\tt 01:} $\mu_{00}(z)\leftarrow 1$;

    \noindent {\tt 02:} {\tt FOR} $\ell = 1$ {\tt TO} $k$ {\tt DO}

    \noindent {\tt 03:} \:\:\: {\tt Compute $\mu_{\ell 0}(z)$ via \eqref{eq:kzero}; }

    \noindent {\tt 04:} {\tt END;}

    \noindent {\tt 05:} {\tt FOR} $n=1$ {\tt TO} $m$ {\tt DO}

    \noindent {\tt 06:} \:\:\:\: {\tt Compute $\mu_{0 n}(z)$ via \eqref{eq:zerom}; }

   \noindent {\tt 07:} \:\:\:\: {\tt FOR} $\ell=1$ {\tt TO} $k$ {\tt DO}

   \noindent {\tt 08:} \:\:\:\:\:\:\:\: {\tt Compute $\mu_{\ell n}(z)$ via \eqref{eq:mukm}; }

   \noindent {\tt 09:} \:\:\:\: {\tt END;}

   \noindent {\tt 10:} {\tt END;}

   \noindent {\tt 11:} {\tt RETURN $\mu_{km}(z)$.}
   
   }
\end{quotation}

\medskip

Observe that in each step of this algorithm, one only needs objects that have been computed earlier.  

\begin{theorem}
    The probability generating function $\mu_{km}(z)$ can be recursively identified via the above algorithm. If $u_{ni}$ and $v_{ni}$, for all $n\in
\{0,\ldots,m\}$ and $i\in\{0,\ldots,n\}$, are known, then the complexity of the algorithm is $O(km^2)$.
\end{theorem}

The complexity of this algorithm is $O(km^2)$ because the complexity of the for-loop in lines ${\tt 07}$--${\tt 09}$ is $O(k n)$, and this loop has to be performed for $n=1$ up to $m$.

\begin{remark}
    One can easily adapt our algorithm to facilitate the computation of the joint transform of the number of customers and the work in the system at the exponential time $T$. Namely, with mild abuse of notations, we first replace $\mu_{\ell n}(z)$ by 
    \begin{equation}
        \mu_{\ell n}(z,\alpha):={\mathbb E}_{\ell n}[z^{Z(T)}{\rm e}^{-\alpha W(T)}],
    \end{equation} where $W(t)$ is the amount of work at time $t$. Note that the quantity $u_{ni}$ remains unchanged, but the $v_{ni}$ needs to be adapted. To this end we observe that, starting from state $(\ell, n)$, with $\ell\in\{1,\ldots,k\}$, then when $T\leqslant B$ and the number of arrivals until time $T$ is $i$, then the total number in the system at time $T$ is $\ell+i$ and the total amount of work at time $T$ is $B-T$ plus a sum of $\ell+i-1$ independent random variables (which are also independent of everything else) that are distributed like the generic service time $B$. Therefore, we see that $z^{\ell+i}v_{ni}$ in \eq{mukm} needs to be replaced by
    \begin{equation}
    z^{\ell+i}\beta(\alpha)^{\ell+i-1}v_{ni}(\alpha)
    \end{equation}
    where $\nu_{00}(\alpha)={\mathbb E}[{\rm e}^{-\alpha (B-T)}1_{\{T\leqslant B\}}]$ and, for $n\in\{1,\ldots,m\}$,
    \begin{equation}
        v_{ni}(\alpha)=
        \begin{cases}
        {\mathbb E}[{\rm e}^{-\alpha(B-T)}1_{\{S_n-S_{n-i}\leqslant T<S_n-S_{n-i-1},T\leqslant B\}}],&i\in\{0,\ldots,n-1\},\\
        {\mathbb E}[{\rm e}^{-\alpha(B-T)}1_{\{S_n\leqslant T,\,T\leqslant B\}}],&i=n.
        \end{cases}
    \end{equation}
As before, we have that $\mu_{00}(z,\alpha)=1$, and \eq{zerom}, \eq{kzero0} and \eq{kzero} remain unchanged, apart from the obvious replacement of $\mu_{\ell n}(z)$ by $\mu_{\ell n}(z,\alpha)$, and of $v_{00}$ by $v_{00}(\alpha)$. In particular, the structure of the algorithm remains the same. 

Clearly, setting $z=1$ (and then $z^{\ell+i}=1$ in the only place where this is relevant), immediately gives an algorithm for the computation of the LST of the workload in the system at an independent exponentially distributed time, starting from the state $(k,m)$.\hfill$\Diamond$
\end{remark}

\begin{remark}
    So far we have assumed that the service times are identically distributed and that the interarrival times $T_m,\ldots,T_1$ are exponentially distributed with parameters $\lambda_m,\ldots,\lambda_1$. This can be substantially generalized, in that we could instead assume that, when after a service completion (or at time zero) there are $n$ customers to arrive and $\ell$ customers present, the next service time is distributed like some $B_{\ell n}$ and that the interarrival times (until the next service completion) have parameters $\lambda_{\ell n},\ldots,\lambda_{\ell 1}$, while maintaining all the underlying independence assumptions.
    
    For this generalization, the only change that is needed in the description is that $u_{ni}$ and $v_{ni}$ are replaced by $u_{kni}$ and $v_{kni}$ where the latter are defined via \eq{umi} and \eq{vmi} in which $B$ is replaced by $B_{kn}$ and $S_i$ are similarly defined only with $T_i\sim {\rm Exp}(\lambda_{ki})$, $i=1,\ldots,m$. All the rest of the development is identical. \hfill$\Diamond$
\end{remark}

\begin{remark}
    We observe that, by differentiation of \eq{mukm}, \eq{zerom} and \eq{kzero0} (or \eq{kzero} directly) $\ell$ times and letting $z\uparrow 1$, then exactly the same algorithm results in the recursive computation of the factorial moments
    \begin{equation}
    \mathbb{E}_{km}(Z(T)(Z(T)-1)\ldots (Z(T)-\ell+1))
    \end{equation}
    for any desired $k,m,\ell$, from which the moments $\mathbb{E}_{km}Z^\ell(T)$ may be deduced in the usual way.

    Similarly, equating the coefficients of the polynomials on the two hand sides of \eq{mukm}, \eq{zerom} and \eq{kzero0}, results in a similar algorithm for the recursive computation of the probabilities $\mathbb{P}_{km}(Z(T)=\ell)$. This is a straightforward exercise. \hfill$\Diamond$
\end{remark}

\begin{remark}
\label{Remark4}
In this remark we study a process that is useful in the context of the waiting-time analysis of Section \ref{sec:wt}. We consider the bivariate Markov chain $(Z_h,V_h)_h$, with $h=1,\ldots,k+m$, where (i)~$Z_h$ denotes the number of customers that is in the system after the $h$-th departure, and (ii)~$V_h$ denotes the number of customers, out of the $m$ that are arriving after time $0$, that still need to arrive after the $h$-th departure. Note that $V_h$ are non-increasing in $h$, almost surely, and that $(Z_0,V_0)=(k,m)$.

In order to analyze the Markov chain $(Z_h,V_h)_h$, we introduce the following sequence of probabilities:
\begin{equation}
w_{ni}:=
\begin{cases}
{\mathbb P}(S_n-S_{n-i}\leqslant B<S_n-S_{n-i-1}),&i\in\{0,\ldots,n-1\}\,,\\
{\mathbb P} (S_n\leqslant B), &i=n\,.
\end{cases}
\end{equation}
Observe that in the limiting regime that $\gamma\downarrow 0$, the probabilities $w_{ni}$ and $u_{ni}$ coincide. 

The transition probabilities of $(Z_h,V_h)_h$ are given by
\begin{equation}
    p_{(\ell_1,n_1),(\ell_2, n_2)}=
    \begin{cases}
    w_{n_1,n_1-n_2},&\ell_1\in\{1,\ldots,k\},\,n_2\in\{0,\ldots,n_1\},\,\ell_2=\ell_1+n_1-n_2-1,\\
    w_{n_1-1,n_1-n_2},&\ell_1=0,\,n_1\in \{1,\ldots,m\},\,n_2\in\{0,\ldots, n_1-1\}, \ell_2=n_1-n_2-1,\\
    {1,}&{\ell_1=n_1=\ell_2=n_2=0},\\
    0&\text{otherwise.}
    \end{cases}
\end{equation}
These cases can be interpreted as follows. The first case corresponds to the scenario that the queue is non-empty at the end of the previous  service time, then a new service time starts, during this service time $n_1-n_2$ customers arrive, with $n_1$ customers that are yet to arrive at the beginning of the service time. 
The second case corresponds to a similar scenario, but now the queue is empty at the end of the previous  service time, and a new service time starts as soon as a new customer enters (so that $n_1-1$ customer arrive during this service time). The third case corresponds to the absorbing state $(0,0)$.

As for $u_{ni}$ and $v_{ni}$ (see Section \ref{sec:spec}), it is easy to check that $w_{ni}$ also have an explicit form for the cases where either $\lambda_i=\lambda$ for $i=1,\ldots,n$ or $\lambda_i=i\lambda$ for $i=1,\ldots,n$. {These transition probabilities will be needed in the next section; more specifically, they are needed to determine the probability that the system is empty after a service completion.}
\hfill$\Diamond$
\end{remark}

\section{Waiting times}\label{sec:wt}
As before, we consider the setting in which $k$ customers are already present at time $0$,
while in addition $m$ customers are yet to arrive after time $0$, with the corresponding interarrival times being represented by $T_m,T_{m-1},\dots,T_1$, where $T_i \sim {\rm Exp}(\lambda_i)$ for $i=1,\dots,m$. We throughout assume a FIFO service discipline, where the initial $k$ customers are arranged at some arbitrary order.
With $W_1,\dots,W_k$ denoting
the waiting times of those customers already present at zero, and measuring their waiting times from time $0$ on, we obviously have 
\begin{equation}
    {\mathbb E}[{\rm e}^{-\alpha W_j}] = \beta^{j-1}(\alpha),\:\:\:j=1,\dots,k.
\end{equation}
The goal of this section is therefore to identify, in terms of transforms, the waiting time distributions of the $m$ customers arriving after time $0$, denoted by $W_{k+1},\ldots,W_{k+m}$.

The waiting times of the customers arriving after time $0$ obey the Lindley recursion 
\begin{equation}W_{j+1} = [W_j+B_j-I_j]^+,\:\:\:
j=k,\dots,k+m-1,  
\end{equation}
where $I_j$ denotes the interarrival time between arriving customers $j-k$ and $j-k+1$
(for $j=k$: the arrival time of customer $k+1$). We can use the identity ${\rm e}^{-\alpha x^+} = {\rm e}^{-\alpha x} - {\rm e}^{-\alpha x^-} + 1$, to conclude that, for $i=0,\dots,m-1$,
\begin{equation}
    {\mathbb E}[{\rm e}^{-\alpha W_{k+i+1}}] =
    {\mathbb E}[{\rm e}^{-\alpha[W_{k+i}+B_{k+i}-I_{k+i}]^+}] = 
    {\mathbb E}[{\rm e}^{-\alpha(W_{k+i}+B_{k+i}-I_{k+i})}]
    - {\mathbb E}[{\rm e}^{-\alpha[W_{k+i}+B_{k+i}-I_{k+i}]^-}] + 1.
    \label{eqW1}
\end{equation}
Using the independence of the components in the first term in the righthand side, and also using the fact that $I_{k+i} \sim {\rm Exp}(\lambda_{m-i})$ (implying that it is memoryless), we obtain from \eqref{eqW1} that
\begin{eqnarray}
    {\mathbb E}[{\rm e}^{-\alpha W_{k+i+1}}] &=&
    {\mathbb E}[{\rm e}^{-\alpha W_{k+i}}] \beta(\alpha) \frac{\lambda_{m-i}}{\lambda_{m-i}-\alpha}
    - {\mathbb P}(W_{k+i+1}=0)\frac{\lambda_{m-i}}{\lambda_{m-i}-\alpha}
-{\mathbb P}(W_{k+i+1}>0)+ 1
\nonumber\\
&=& {\mathbb E}[{\rm e}^{-\alpha W_{k+i}}] \beta(\alpha) \frac{\lambda_{m-i}}{\lambda_{m-i} - \alpha} - {\mathbb P}(W_{k+i+1}=0) \frac{\alpha}{\lambda_{m-i}-\alpha} .
\label{eqW2}
\end{eqnarray}
Iterating this relation, we readily obtain the following result. In the sequel we use the notation $\varrho_h:= {\mathbb P}(W_h=0)$, for $h=k+1,\ldots,k+m.$

\begin{theorem}
For $j=k+1,\dots,k+m$, empty products being equal to one,
\begin{equation}
    {\mathbb E}[{\rm e}^{-\alpha W_j}] = 
    \beta^{j-1}(\alpha) \prod_{i=0}^{j-k-1} \frac{\lambda_{m-i}}{\lambda_{m-i}-\alpha}
    - \sum_{h=k+1}^{j
    }  \varrho_h \beta^{j-h}(\alpha) \frac{\alpha}{\lambda_{m-h+k+1} - \alpha} 
    \prod_{w=0}^{j-1-h} \frac{\lambda_{m+k-h-w}}{\lambda_{m+k-h-w}-\alpha} .
    \label{eqW3}
\end{equation}
\end{theorem}
We have thus derived an expression for the waiting time LST of all $k+m$ customers,
in known quantities and in the constants $\varrho_h$ for $h=k+1,\dots,k+m$.
Now observe that $\varrho_h = {\mathbb P}(Z_{h-1}=0)$.
Hence one only needs to determine the latter probabilities from the one-step transition probabilities ${\mathbb P}(Z_{i+1}=j\,|\,Z_i=r)$ which are discussed in Remark~\ref{Remark4} above.

Differentiating (\ref{eqW3}) with respect to $\alpha$ and subsequently taking $\alpha=0$, we obtain that if ${\mathbb E}[B]<\infty$, then, for $j=k+1,\dots,k+m$:
\begin{equation}
    {\mathbb E}[W_j] = (j-1) {\mathbb E}[B] - \sum_{i=0}^{j-k-1} \frac{1}{\lambda_{m-i}} + \sum_{h=k+1}^j \varrho_h \frac{1}{\lambda_{m-h+k+1}}.
    \label{eqW4}
    \end{equation}
Conclude in particular that if ${\mathbb E}[B]<\infty$, then $ {\mathbb E}[W_j]<\infty$.

We end this section by exploiting (\ref{eqW3}) and (\ref{eqW4})
to derive the tail asymptotics of the distribution of waiting time $W_j$ in the case of a regularly varying service-time distribution, an important class of heavy-tailed distributions. We use the well-known Bingham-Doney lemma, cf.\ \cite[Theorem 8.1.6]{BGT}, which states that (i) and (ii) below are equivalent: (i) $B$ is regularly varying of index $-\nu \in (-2,-1)$ at infinity, i.e.,
\begin{equation}
{\mathbb P}(B>t) \sim \frac{-1}{\Gamma(1-\nu)} t^{-\nu} L(t), ~~~ t \rightarrow \infty,
\label{eqW5}
\end{equation}
with $L(t)$ a slowly-varying function (entailing that the first moment of $B$ is finite but the second is not), and (ii)
\begin{equation}
\beta(\alpha) -1 + \alpha {\mathbb E}[B] \sim \alpha^{\nu} L(1/\alpha), ~~~ \alpha \downarrow 0. 
\label{eqW6}
\end{equation}
Now assume that (\ref{eqW5}) and hence also (\ref{eqW6}) hold.
It then readily follows from (\ref{eqW3}) and (\ref{eqW4}) that 
\begin{equation}
{\mathbb E}[{\rm e}^{-\alpha W_j}] -1 + \alpha {\mathbb E}[W_j] \sim (j-1) \alpha^{\nu} L(1/\alpha), ~~~ \alpha \downarrow 0.
\label{eqW7}
\end{equation}
Once more applying the equivalence between tail behaviour of the distribution and behavior of the LST near zero that is specified by the Bingham-Doney lemma, but now for $W_j$,
we find the following result. We denote by $f(t)\sim g(t)$ that $f(t)/g(t)\to 1$ as $t\to\infty.$

\begin{proposition}
    Suppose the service-time distribution obeys (\ref{eqW5}). Then, 
\begin{equation} 
    {\mathbb P}(W_j>t) \sim (j-1) {\mathbb P}(B>t), ~~~ t \rightarrow \infty .
    \label{eqW9}
    \end{equation}

\end{proposition}
    For the ordinary M/G/1 queue in steady state, with regularly varying service time distribution as in (\ref{eqW5}), it is well-known \cite{Cohen73} that the waiting time is also regularly varying but with index $1-\nu$.
    However, in our setting with only a finite number of customers, the most likely scenario for a very large waiting time $W_j$ is that one of the preceding $j-1$ service times is very large, resulting in (\ref{eqW9}); 
    in this respect, bear in mind that for independent $B_1,\ldots, B_j$ that are distributed as the random variable $B$, we have that
    \begin{equation}
        {\mathbb P}\left(\sum_{i=1}^j B_i > t\right) \sim j\,{\mathbb P}(B>t), ~~~ t\to\infty.
    \end{equation}

\section{Special cases}  \label{sec:spec}
We now consider two special cases in which the required probability sequences $u_{mi}$ and $v_{mi}$ allow an explicit evaluation. 
In the first case we take $\lambda_i=i\lambda$. This means that the random vector $(T_m,\ldots,T_1)$ is distributed as the differences of the components of the order statistics associated with $m$ i.i.d.\ ${\rm Exp}(\lambda)$ distributed random variables. This arrival process is the pure death process analogue of the Yule process (which is a pure birth process). We can treat it as the case where the $m$ customers who are yet to arrive pick their arrival times according to independent (identically distributed) exponential clocks. Whereas for the Poisson process the inter-arrival times are i.i.d.\  exponentially distributed, for the case with $\lambda_i=i\lambda$ the arrival times are i.i.d.\ exponentially distributed.
Hence, as is easily checked, for $i\in\{0,\ldots, m\}$,
\begin{equation}
u_{mi}=
\binom{m}{i}\,{\mathbb E}\left({\rm e}^{-\gamma B}(1-{\rm e}^{-\lambda B})^i\,{\rm e}^{-(m-i)\lambda B}\right)\,,
\end{equation}
To compute $v_{mi}$ we again recall \eq{vmi}, \eq{vmi1}, \eq{vmi2}, \eq{vmi3} and \eq{vmi4}, so as to obtain, for $i\in\{0,\ldots,m\}$,
\begin{equation}
v_{mi}=\frac{\gamma}{\gamma +(m-i)\lambda}\,\prod_{j=m-i+1}^{m}\frac{j\lambda}{\gamma +j\lambda}-
\binom{m}{i}\,
{\mathbb E}\left({\rm e}^{-\gamma B}(1-{\rm e}^{-\lambda(B+T)})^i\, {\rm e}^{-(m-i)\lambda (B+T)}\right)\,.
\end{equation}
The expressions for $u_{mi}$ and $v_{mi}$ can be made more explicit by noting that
\begin{align}
    {\mathbb E}\left({\rm e}^{-\gamma B}(1-{\rm e}^{-\lambda B})^i\,{\rm e}^{-(m-i)\lambda B}\right)&=\sum_{j=0}^i\binom{i}{j}(-1)^{i-j}\beta(\gamma+(m-j)\lambda),\\
    {\mathbb E}{\rm e}^{-\gamma B}\left((1-{\rm e}^{-\lambda(B+T)})^i\, {\rm e}^{-(m-i)\lambda (B+T)}\right)&=\sum_{j=0}^i\binom{i}{j}(-1)^{i-j}\beta(\gamma+(m-j)\lambda)
    \,\frac{\gamma}{\gamma+(m-j)\lambda}\,.
\end{align}

\medskip

    The remainder of this section covers a second case in which $\lambda_i\equiv \lambda$, i.e., arrivals occur according to a Poisson process with rate $\lambda$ which is stopped after the $m$-th arrival. By Eqn.\ \eq{umi}, 
\begin{equation}\label{eq:poissonu}
u_{mi}=
\begin{cases}
{\mathbb E}\left({\rm e}^{-(\lambda+\gamma )B}{\displaystyle \frac{(\lambda B)^i}{i!}}\right),&i\in\{0,\ldots, m-1\},\\
\beta(\gamma)-\sum_{j=0}^{m-1}u_{mj},&i=m\,.
\end{cases}
\end{equation}
Also, recalling \eq{vmi}, \eq{vmi1}, \eq{vmi2}, \eq{vmi3} and \eq{vmi4}, 
\begin{equation}\label{eq:poissonv}
v_{mi}=
\begin{cases}
\left({\displaystyle \frac{\lambda}{\lambda+\gamma }}\right)^i\,{\displaystyle \frac{\gamma }{\gamma +\lambda}}-{\mathbb E}\left({\rm e}^{-((\lambda+\gamma)B+\lambda T)}{\displaystyle \frac{(\lambda (B+T))^i}{i!}}\right),&i\in\{0,\ldots, m-1\}\\
1-\beta(\gamma)-\sum_{j=0}^{m-1}v_{mj},&i=m\,.
\end{cases}
\end{equation}
Observe that in \eq{poissonu} and \eq{poissonv} we can also write, for $i\in\{0,\ldots,m-1\}$,
\begin{align}
    \mathbb{E}\left({\rm e}^{-(\lambda+\gamma)B}\frac{(\lambda B)^i}{i!}\right)&=\frac{(-\lambda)^i\beta^{(i)}(\lambda+\gamma)}{i!}\\
    {\mathbb E}\left({\rm e}^{-((\lambda+\gamma)B+\lambda T)}{\displaystyle \frac{(\lambda (B+T))^i}{i!}}\right)&=\sum_{j=0}^i{\mathbb E}\left({\rm e}^{-(\lambda+\gamma) B}{\displaystyle \frac{(\lambda B)^j}{j!}}\right){\mathbb E}\left({\rm e}^{-\lambda T}{\displaystyle \frac{(\lambda T)^{i-j}}{(i-j)!}}\right)
\notag \\&
    =\sum_{j=0}^i \frac{(-\lambda)^j\beta^{(j)}(\lambda+\gamma)}{j!}\,\frac{\gamma}{\lambda+\gamma}\left(\frac{\lambda}{\lambda+\gamma}\right)^{i-j}\,.
\end{align}

Throughout this paper we have assumed that $k$ customers are initially present, and there is a pool of $m$ customers who arrive after time zero. Below we briefly study the case in which both $k$ and $m$ are, geometrically distributed, random variables.
Restricting ourselves to the case $\lambda_i \equiv \lambda$ and to Exp($\mu$) distributed service times, we obtain an explicit expression for the generating function of the number of customers at an Exp($\gamma$) distributed time (when divided by $\gamma$, this is also the Laplace transform of the generating function of the number of customers at time $t$). This generalizes a known result (cf.\ \cite[Ch.\ II.2]{Cohen82}) for the transient behavior of the ordinary M/M/1 queue to the case in which the arrival process dries up after a geometrically distributed number of arrivals.
In this case, recalling Remark~\ref{remBexp} and with $\xi\equiv\xi(\lambda,\mu,\gamma):=\lambda+\mu+\gamma$, \eq{poissonu} and \eq{poissonv} become

\begin{align}
u_{mi} &=
\begin{cases}
     \frac{\mu}{\LMG} \left(\frac{\lambda}{\LMG}\right)^i, &i\in\{0,1,\ldots,m-1\},\\
     \frac{\mu}{\MG} \left(\frac{\lambda}{\LMG}\right)^m, & i=m\,,
    \end{cases}
    \label{a.1}
\end{align}
and
\begin{align}
v_{mi} &=
\begin{cases}
     \frac{\gamma}{\LMG} \left(\frac{\lambda}{\LMG}\right)^i, &i\in\{0,1,\ldots,m-1\},\\
      \frac{\gamma}{\MG} \left(\frac{\lambda}{\LMG}\right)^m, & i=m\,.
    \end{cases}
    \label{a.2}
\end{align}
The above equations have obvious interpretations. For instance, regarding \eqref{a.1}: among three Poisson events with rates $\lambda,\mu,\gamma$, an arrival occurs $i$ times first with probability $(\lambda/\xi)^ i$, and  a service completion comes next with probability $\mu/\xi$
(but with probability $\mu/(\mu+\gamma)$ if there is no potential arrival left).

Introduce, for $p,r\in[0,1]$ and $z\in(0,1]$,
\begin{align}
G(p,r,z) &:= \sum_{\ell=0}^{\infty} p^\ell  \sum_{n=0}^{\infty} r^n \mu_{\ell,n}(z),\quad\quad
    M_\ell(r,z) := \sum_{n=0}^{\infty} r^n \mu_{\ell,n}(z),\quad \ell=0,1. 
    \label{a.3}
    \end{align}
    Observe that $G(p,r,z)$ is the generating function of $\mu_{\ell,n}(z)$ in both $\ell$ and $n$, but it should be multiplied by $(1-p)(1-r)$ to obtain the generating function of $Z(T)$  for geometrically distributed $k$ and $m$.  {Note that, applying numerical inversion to $G(p,r,z)$, we can identify $\mu_{\ell,n}(z)$, and hence we can in principle also deal with arbitrarily distributed $k$ and $m$.}
    In the rest of this section we find a closed-form expression for $G(p,r,z)$.
Taking generating functions in (\ref{eq:mukm}) yields:
\begin{equation}
G(p,r,z) - M_0(r,z) = 
\sum_{\ell=1}^{\infty} p^\ell  \sum_{n=0}^{\infty} r^n \sum_{i=0}^n \mu_{\ell+i-1,n-i}(z) \,u_{ni} +
\sum_{\ell=1}^{\infty} p^\ell  \sum_{n=0}^{\infty} r^n \sum_{i=0}^n z^{\ell+i}\, v_{ni} =: {\rm I + II}.
\label{a.5}
\end{equation}
Recall the identity
$(a-b)\sum_{j=0}^{\infty} a^j \sum_{i=0}^{\infty} b^i c_{j+i} = {a \sum_{h=0}^{\infty} c_h a^h - b \sum_{h=0}^{\infty} c_h b^h}$ for any $a,b\in[0,1]$ and a sequence $c_{k}$. Hence, 
\begin{align}
   G(p,r,z) - M_0(r,z) =&\: \frac{\mu p}{\LMG - \lambda r/p} G(p,r,z) -
   \frac{\lambda r}{\LMG} \frac{\mu}{\LMG - \lambda r/p} G\left(\frac{\lambda r}{\LMG},r,z\right)\:+
   \nonumber 
   \\
   &\frac{p \mu}{\MG-p\mu} H_1(r,z)
   - \frac{pz}{1-pz} H_2(r,z)+ \frac{pz}{1-pz} 
   F(r,z),
    \label{a.7}
\end{align}
where
\begin{align}
H_1(r,z)&:= \frac{(\MG)(1-z)}{\mu - (\MG)z}
\frac{\lambda \mu}{(\MG) \LMG -r\mu\lambda},
\label{a.8}\\
    H_2(r,z)& := \frac{\gamma}{\mu-(\MG)z} \frac{\mu}{\MG} \frac{\lambda}{\LMG-\lambda rz},
    \label{a.9}\\
    F(r,z)&:=
 \frac{\gamma}{\MG} \frac{1}{1-r} \frac{\xi-\lambda r}{\xi -\lambda rz}. 
   \label{a.10}
   \end{align}
Here the first four terms in the right-hand side of (\ref{a.7}) concern I, and the fifth term II.  The third and fourth term correct for taking the expression in the first (instead of the second) line of \eqref{a.1} for $u_{nn}$.
Observe that $z = \mu/(\mu+\gamma)$ is a pole of both the third and the fourth term, but {\em not} of their difference.
In the fifth term we have used both parts of \eqref{a.2}, and also the first line of (\ref{eq:kzero}). When $z=u_{00}$, hence $z = {\mu}/{\xi}$, one should use the second line of (\ref{eq:kzero});
we have omitted that calculation.

The next step is to rewrite (\ref{a.7}) into
\begin{align}
    G(p,r,z) \left(\frac{\mu p^2 - \LMG p +\lambda r}{\lambda r -  \LMG p}\right)& = M_0(r,z) + \frac{\lambda r}{\LMG} \frac{\mu p}{\lambda r -  \LMG p} G\left(\frac{\lambda r}{\LMG},r,z\right)
    \nonumber
    \\
    &\:\:+\: \frac{p\mu}{\MG-p \mu} H_1(r,z) - \frac{pz}{1-pz} H_2(r,z) + \frac{pz}{1-pz} F(r,z).
    \label{a.11}
\end{align}
It is readily verified that the numerator of the term between brackets in the left-hand side of (\ref{a.11}) has one zero, say $p^\star $, in $(0,1)$ (use Rouch\'e's theorem, or take $p=0$ and $p=1$ in that numerator), and a second zero $\hat{p} = {\lambda r}/({\mu}{p^\star }) > 1$.
Since $G(p^\star ,r,z)$ is finite, the right-hand side of (\ref{a.11}) must be zero for $p=p^\star $, yielding a relation between the unknown functions $G(\frac{\lambda r}{\LMG},r,z)$ and $M_0(r,z)$:
with $J(r,z):= H_2(r,z) - F(r,z)$,
\begin{align}
    \frac{\lambda r}{\LMG} \frac{\mu p^\star }{\lambda r -  \LMG p^\star } G\left(\frac{\lambda r}{\LMG},r,z\right) =\:&
   -M_0(r,z) - \frac{p^\star  \mu}{\MG - p^\star  \mu} H_1(r,z) 
   +\frac{p^\star  z}{1-p^\star  z} J(r,z).
   \label{a.12}
   \end{align}
   Substituting this in (\ref{a.11}) gives, after some calculations,
   \begin{align}
       \mu (p-\hat{p}) G(p,r,z) =&\: - \frac{\lambda r}{p^\star } M_0(r,z) 
       + \frac{p\mu}{\MG-p\mu} \frac{\lambda \mu r - (\MG) \LMG }{(\MG-p^\star \mu)} H_1(r,z)\:-\notag\\
       &\:\frac{pz}{1-pz}  \frac{\lambda r z -  \LMG }{1-p^\star z}J(r,z) .
       \label{a.13}
   \end{align}
   It remains to determine $M_0(r,z)$. Using the one equation we have not used yet,  namely (\ref{eq:zerom}), in combination with $\mu_{0,0}(z) \equiv 1$, we obtain
   \begin{equation}
       M_0(r,z) -1 = \frac{\gamma}{\LG} \frac{r}{1-r} + \frac{\lambda}{\LG} r M_1(r,z).
       \label{a.14}
       \end{equation}
       From (\ref{a.3}), in combination with (\ref{a.14}), we see that
       \begin{equation}
           \lim_{p \downarrow 0} \frac{G(p,r,z)-M_0(r,z)}{p} = 
           M_1(r,z) = \frac{\LG}{\lambda r} \left(M_0(r,z) -1-\frac{\gamma}{\LG} \frac{r}{1-r}\right).
           \label{a.15}
           \end{equation}
        Subtracting $\mu(p-\hat{p}) M_0(r,z)$ from both sides in (\ref{a.13}), dividing by $p$ and finally letting $p \downarrow 0$ results, in combination with (\ref{a.15}), in an equation from which $M_0(r,z)$ can be determined  
         (notice that, in combining the resulting two $M_0(r,z)$ terms in the right-hand side, we have used that $p^\star  \hat{p} = \lambda r/\mu$):
           \begin{align}
                \frac{\LG}{\lambda r}& \left(M_0(r,z)-1-\frac{\gamma}{\LG} \frac{r}{1-r}\right)\nonumber
            \\
            &=\frac{M_0(r,z)}{\hat{p}} -\frac{1}{\mu \hat{p}} \left(\frac{\mu}{\MG} \frac{\lambda \mu r - (\MG) \LMG }{\MG - p^\star \mu} H_1(r,z) -z \frac{\lambda r z -  \LMG }{1-p^\star z} J(r,z)\right).
             \label{a.16}
            \end{align}
            Hence
            \begin{align}
                 M_0(r,z) = &\:\frac{(\LG) \hat{p}}{(\LG) \hat{p} - \lambda r} + \frac{\gamma r \hat{p}}{(1-r)((\LG) \hat{p} - \lambda r}\:-
                \nonumber\\
                &\: \frac{\lambda r}{\mu} \frac{1}{(\LG) \hat{p} - \lambda r} 
            \left(\frac{\mu}{\MG} \frac{\lambda \mu r - (\MG) \LMG }{\MG - p^\star \mu} H_1(r,z) -z \frac{\lambda r z -  \LMG }{1-p^\star z} J(r,z)\right).  
            \label{a.17}
            \end{align}
            Notice that $M_0(0,z) = 1$, as it should be.

\smallskip
            
           In conclusion, we find the following result. Recall that $\xi\equiv\xi(\lambda,\mu,\gamma)=\lambda+\mu+\gamma$.
    \begin{proposition} For any $p,r\in[0,1]$ and $z\in(0,1]$, $G(p,r,z)$ is given by \eqref{a.13}, with $M_0(r,z)$ given by \eqref{a.17}.
    \end{proposition}

\end{document}